\documentclass{article}
\usepackage[comma,numbers,square,sort&compress]{natbib}

\usepackage{subfigure}
\usepackage{graphicx,amsmath,amssymb}
\usepackage{subfigure,pgf,tikz,pst-node}
\usetikzlibrary{arrows,shapes,automata,positioning,fit}
\usepackage{multirow}
\newtheorem{theorem}{Theorem}
\newtheorem{lemma}{Lemma}

\newtheorem{assumption}{Assumption}
\newtheorem{remark}{Remark}

\begin{document}

\title{Coordination of Multi-Agent Systems under Switching Topologies via Disturbance Observer Based Approach}

\author{Yutao Tang $^{\ast}$\thanks{$^\ast$Yutao Tang is with School of Automation, Beijing University of Posts and Telecommunications, Beijing 100876, China.  E-mail\,$: yttang@amss.ac.cn$}}
\date{}
\maketitle

\begin{abstract}
In this paper, a leader-following coordination problem of heterogeneous multi-agent systems is considered under switching topologies where each agent is subject to some local (unbounded) disturbances. While these unknown disturbances may disrupt the performance of agents, a disturbance observer based approach is employed to estimate and reject them. Varying communication topologies are also taken into consideration, and their byproduct difficulties are overcome by using common Lyapunov function techniques. According to the available information in difference cases, two disturbance observer based protocols are proposed to solve this problem. Their effectiveness is verified by simulations.


\end{abstract}

\section{Introduction}

In the past decades there has been a large literature in the study of multi-agent systems due to its wide applications such as cooperative control of unmanned aerial vehicles, communication among sensor networks, and formation of mobile robots (see \cite{olfati2007consensus,ren2008distributed} and the references therein). As an important topic of multi-agent systems, the leader-following problem is actively studied by many authors, e.g. \cite{leonard2001virtual, ren2005consensus, clark2014minimizing}. In this formulation, one or multiple agents are selected as leaders to generate desired trajectories for those followers and lead the whole group to achieve collective tasks \cite{hong2006tracking, zhang2011optimal, lou2012target, meng2013global, peng2014distributed, wen2013consensus, tang2015robust}.

It has been well-recognized that disturbance rejection is of fundamental importance in the applicability of designed controllers. While there are always disturbances in real applications, it is necessary to take them into consideration and attenuate or eliminate them for multi-agent control design. In \cite{bauso2009consensus} for first-order multi-agent systems, it was shown that under bounded unknown external disturbances the steady-state errors of any two agents can reach a small region and is called lazy consensus. Later, the authors of \cite{liu2010consensus} proposed an $H_\infty$ analysis approach to investigate robust consensus problem of high-order multi-agent systems with external disturbances. In \cite{zhang2013further}, consensus of multi-agent systems with exogenous disturbances was considered, and an observer was constructed to compensate the negative effect of those disturbances. However, most of these results were obtained for special dynamic systems, and there are few general consensus results emphasizing disturbance rejection with an exception \cite{su2012scl}, where the authors considered leader-following consensus with disturbance rejection from the viewpoint of output regulation \cite{wonham1979linear} and solved it for linear multi-agent systems with a fixed topology.

Also note that, in centralized or decentralized setup, there is usually no need to partition and treat those disturbances and references in a separate way. However, in multi-agent systems, it would be better to distinguish those two kinds of signals and not to model them in the same manner. Intuitively, the reference is globally set up to drive all agents to complete a common task, while the disturbances are usually local and harmful to such cooperation. Thus, another motivation of this paper is to formulate a problem that treats those two kinds of signals by different approaches.

Hence, we aim to investigative a coordination problem of heterogeneous multi-agent systems under switching topologies, where the reference is given by a traditional leader. At the same time, those followers may be subject to local disturbances modeled by some other autonomous systems. Since the disturbances are often unmeasurable, a disturbance observer based (DOB) approach \cite{nakao1987robust, chen2000nonlinear} is employed to tackle this problem. {DOB approach stems from feedforward control, and can be perceived as a composite controller comprising a feedforward compensation part to reject those disturbances, based on disturbance observation and a feedback rule to regulate the plant to achieve other goals. Although it has been investigated by many publications (see \cite{li2014disturbance} and references therein), there is no corresponding result to our knowledge for heterogeneous multi-agent systems.}

To sum up, the main contributions of the present paper are at least twofold:
\begin{itemize}
\item We extend the conventional leader-following consensus \cite{hong2006tracking, ren2008distributed} to general linear multi-agent systems with local disturbances (which may be unbounded). When there are no such disturbances, these results are consistent with existing consensus results.   Here we consider heterogeneous multi-agent systems under switching topologies, while many existing results were derived for only second order dynamic systems \cite{hu2010distributed} or for fixed graph cases \cite{su2012scl}.
\item We extend the conventional disturbance observer based (DOB) approach \cite{chen2000nonlinear} to its distributed version for multi-agent systems with both reference tracking and disturbance rejection. When there is only one agent, our problem becomes the conventional DOB formulation. Even for the centralized case, we propose different full-order and reduced-order disturbance observers to solve this problem without using the derivative of the plant's states as that in \cite{li2014disturbance}. It is also remarkable that these disturbance observers can allow both bounded disturbances (e.g., constant and harmonic signals in existing literature) and unbounded disturbances (e.g. ramping and polynomial signals).
\end{itemize}
The rest of this paper is organized as follows. In Section 2, some preliminaries are given and our problem is formulated. Then main results are presented in Section 3, where two types of control laws are constructed. Finally, simulations and our concluding remarks are provided in Sections 4 and 5, respectively.

Notations:   Let $\mathbb{R}^n$ be the $n$-dimensional Euclidean space, $\mathbb{R}^{n\times m}$ be the set of $n\times m$ real matrices. $\mathbf{0}^{n \times m}$ represents an $n\times m$ zero matrix. $\text{diag}\{b_1,{\dots},b_n\}$ denotes an $n\times n$ diagonal matrix with diagonal elements $b_1,\,\cdots\,,b_n$; $\text{col}(a_1,{\dots},a_n) = [a_1^T,{\dots},a_n^T]^T$ for column vectors $a_i\; (i=1,{\dots},n)$. A weighted directed graph (or weighted digraph) $\mathcal {G}=(\mathcal {N}, \mathcal {E}, \mathcal{A})$ is defined as follows, where $\mathcal{N}=\{1,{\dots},n\}$ is the set of nodes, $\mathcal {E}\subset \mathcal{N}\times \mathcal{N}$ is the set of edges, and $\mathcal{A}\in \mathbb{R}^{n\times n}$ is a weighted adjacency matrix \cite{mesbahi2010graph}. $(i,j)\in \mathcal{E}$ denotes an edge leaving from node $i$ and entering node $j$. The weighted adjacency matrix of this digraph $\mathcal {G}$ is described by $A=[a_{ij}]_{i,\,j=1,\dots,n}$, where $a_{ii}=0$ and $a_{ij}\geq 0$ ($a_{ij}>0$ if and only if there is an edge from agent $j$ to agent $i$). A path in graph $\mathcal {G}$ is an alternating sequence $i_{1}e_{1}i_{2}e_{2}{\cdots}e_{k-1}i_{k}$ of nodes $i_{l}$ and edges $e_{m}=(i_{m},i_{m+1})\in\mathcal {E}$ for $l=1,2,{\dots},k$. If there exists a path from node $i$ to node $j$ then node $i$ is said to be reachable from node $j$. The neighbor set of agent $i$ is defined as $\mathcal{N}_i=\{j\colon (j,i)\in \mathcal {E} \}$ for $i=1,...,n$.  A graph is said to be undirected if $a_{ij}=a_{ji}$ ($i,j=1,{\dots},n$). The weighted Laplacian $L=[l_{ij}]\in \mathbb{R}^{n\times n}$ of graph $\mathcal{G}$ is defined as $l_{ii}=\sum_{j\neq i}a_{ij}$ and $l_{ij}=-a_{ij} (j\neq i)$.

\section{Problem formulation}

In this paper, we consider $N+1$ agents and $N$ of them are followers of the form:
\begin{align}\label{follower}
\begin{split}
  \dot{x}_i&=A_ix_i+B_iu_i+E_id_i\\
  y_i&=C_ix_i+D_iu_i, \quad i=1,\dots, N
\end{split}
\end{align}
where $x_i\in \mathbb{R}^{n_i}$, $y_i\in \mathbb{R}^{l}$, and $u_i\in \mathbb{R}^{m_i}$ are the state, output, and input of the $i$th subsystem, respectively.  $d_i\in \mathbb{R}^{q_i}$ is the local disturbance of agent $i$ governed by
\begin{align}\label{disturbance}
  \dot{d}_i=S_id_i.
\end{align}
The reference signal is given by a leader (denoted as agent $0$) described as
\begin{align}\label{leader}
  \dot{r}=S_0r, \quad y_0=F_0r,\quad r\in \mathbb{R}^{n_0}
\end{align}
Without loss of generality, we assume $(F_0, S_0)$ is detectable and $S_0$,\,\dots,\,$S_N$ have no eigenvalues lying in the open left half plane.  Let $e_i=y_i-y_0 \,(i=1, \dots, N)$, we aim to design proper controllers such that for any initial condition $x_i(0),\,d_i(0), r(0)$, the tracking error $e_i$ will converge to zero as time goes to infinity in spite of these disturbances.

Unlike in centralized cases, we do not assume the availability of $y_0$ (hence $r$) to all agents in our problem. An agent can get access to $y_0$ unless there is an edge between this agent and the leader. This makes it much difficult to achieve collective behaviors. Associated with this multi-agent system, a dynamic digraph $\mathcal{G}$ can be defined with the nodes $\mathcal{N}=\{0,1,..., N\}$ to describe the communication topology, which may be switching. If the control $u_i$ can get access to the information of agent $j$ at time instant $t$, there is an edge $(j,i)$ in the graph $\mathcal{G}$, i.e., $a_{ij}>0$.  Also note that $a_{0i}=0$ for $i=1,...,N$, since the leader won't receive any information from the followers. Denote the induced subgraph associated with all followers as $\bar{\mathcal{G}}$.

We say a communication graph is connected \cite{hong2006tracking} if the leader (node 0) is reachable from any other node of $\mathcal{G}$ and the induced subgraph $\bar{\mathcal{G}}$ is undirected. Given a communication graph $\mathcal{G}$, denote $H\in \mathbb{R}^{N\times N}$ as the submatrix of its Laplacian $L$ by deleting the first row and first column. By Lemma 3 in \cite{hong2006tracking}, $H$ is positive definite if the communication graph is connected. Denote its eigenvalues as $\lambda_1\geq\lambda_2\geq{\cdots}\geq \lambda_N>0$.

In multi-agent systems, the connectivity graph $\mathcal{G}$ may be time-varying. To describe the variable interconnection topology, we denote all possible communication graphs as $\mathcal{G}_1$,\ldots,$\mathcal{G}_\kappa$, $\mathcal{P}=\{1,\dots, \kappa\}$, and define a switching signal $\sigma: [0,\infty)\rightarrow \mathcal{P}$, which is piece-wise constant defined on an infinite sequence of nonempty, bounded, and contiguous time-intervals. Assume  $t_{i+1}-t_i\geq \tau_0>0, \;\forall i$, where $t_{i}$ is the $i$th switching instant and $t_0=0$. Here $\tau_0$ is often called the dwell-time. Therefore, $\mathcal{N}_i$ and the connection weight $a_{ij}\; (i,j=0,1,\ldots,N)$ are time-varying. Moreover,  the Laplacian $L_{\sigma(t)}$ associated with the switching interconnection graph $\mathcal{G}_{\sigma(t)}$ is also time-varying (switched at $t_i,\; i=0,1,\ldots$), though it is a time-invariant matrix in each interval $[t_{i},t_{i+1})$. The following assumption on the communication graph is often made \cite{hong2006tracking}.
\begin{assumption}\label{ass:graph}
The graph $\mathcal{G}_{\sigma(t)}$ is switching among a group of connected graphs.
\end{assumption}

The coordination problem of these heterogeneous multi-agent systems composed of \eqref{follower}, \eqref{disturbance}, and \eqref{leader} under switching topologies is described as follows. {\it Given these multi-agent systems and the communication graph $\mathcal{G}_{\sigma(t)}$, find a proper distributed control law such that for all initial conditions of the closed-loop system, we have \begin{align}\label{eq:regulate-goal} \lim\limits_{t\to+\infty}e_i(t)=0, \quad i=1,\dots, N.\end{align}}

\begin{remark}\label{rem:dor}
 While a large literature in multi-agent systems only considered a globally reference tracking problem \cite{ni2010leader,tang2014leader} or treated those disturbances and the referenced signal in the same manner \cite{su2012scl}, both reference tracking and disturbance rejection are considered in this formulation while local disturbances are modeled by separate autonomous systems. When $d_i=0$, this formulation is consistent with existing consensus results for general linear dynamics, e.g., \cite{ni2010leader}, which include the well-known consensus for integrators \cite{olfati2007consensus} as its special case.
\end{remark}

\section{Main results}

In this section, we employ a two-phase design procedure to achieve the coordination goal. First, we construct a distributed observer for each agent and transform such a coordination problem of this multi-agent system into $N$ decentralized sub-tasks. Then, those sub-tasks will be completed via DOB approach together with both full-order and reduced-order controllers.

The following lemma is useful and can be proved from the convergent-input convergent-state property of the first subsystem \cite{sontag2008input}.
\begin{lemma}\label{lem:time-varying}
  Consider the linear time-varying system
  \begin{eqnarray*}
  \begin{split}
    \dot{\bar x}&=\bar{A} \bar x+\bar B \bar v\\
    \dot{\bar v}&=\bar S(t)\bar v\\
    \bar e&=\bar{C} \bar{x}+\bar{F}\bar v.
  \end{split}
  \end{eqnarray*}
If $\bar A$ is Hurwitz and the subsystem $\dot{\bar v}=\bar S(t)\bar v$ is uniformly exponentially stable, then for any  $\bar{x}(0)=\bar{x}_0$ and $\bar v(0)=\bar v_0$, it holds $ \lim_{t\to+\infty}\bar e(t)=0$.
\end{lemma}

By letting $v_i\triangleq\mbox{col}(d_i, r)$ and some mathematical manipulations, this leader-following coordination problem is equivalent to the following $N$ sub-problems:
\begin{align}\label{sys:virtual-regulation}
\begin{split}
  \dot{x}_i&=A_ix_i+B_iu_i+\bar E_iv_i\\
  \dot{v}_i&=\bar S_i v_i\\
  e_i&=C_ix_i+D_iu_i+\bar F_iv_i
  \end{split}
\end{align}
where $\bar E_i=[E_i, \,\mathbf{0}^{n_i\times l}]$, $\bar S_i=\mbox{block}\,{diag}\{S_i,\,S_0\}$, and $\bar F_i=[\mathbf{0}^{l\times q_i},\,-F_0]$. $v_i$ is the exogenous signal for agent $i$ including its local disturbance $d_i$ and the global reference $r$.

The following equations known as regulator equations \cite{huang2004nonlinear} play a key role in solving the coordination problem of multi-agent systems.
\begin{assumption}\label{ass:regulator-equation}
  For each $i=1,\,\dots,\,N$, there exist constant matrices $X_{i1}$, $X_{i2}$, $U_{i1}$, and $U_{i2}$ satisfying
  \begin{align*}
    \begin{split}
      X_{i1}S_i&=A_iX_{i1}+B_iU_{i1}+E_i\\
      0&=C_iX_{i1}+D_iU_{i1}
    \end{split}
    \end{align*}
    and
    \begin{align*}
    \begin{split}
      X_{i2}S_0&=A_iX_{i2}+B_iU_{i2}\\
      0&=C_iX_{i2}+D_iU_{i2}-F_0.
    \end{split}
  \end{align*}
\end{assumption}
\begin{remark}\label{rem:transmission-zero}
  Similar conditions have been used in \cite{wang2013distributed,su2012scl}. A sufficient condition to the solvability of these linear matrix equations is that, for any eigenvalue of $S_i$ (denoted as $\lambda$), $i=0,\dots,N$,
  \begin{align*}
      \text{rank}\begin{bmatrix}
        A_i-\lambda I&    B_i\\
        C_i          &    D_i
      \end{bmatrix}=n_i+p.
  \end{align*}
  Specially, when all agents are homogenous without local disturbances, $X_{i2}$ can be taken as the identity matrix which was implicitly used in \cite{ren2008distributed} and \cite{meng2013global}.
\end{remark}

Note that $(A_i, B_i)$ is stabilizable, there exists $K_{i1}$ such that $A_i+B_iK_{i1}$ is Hurwitz. Denote $K_{i2}\triangleq U_{i1}-K_{i1}X_{i1},\, K_{i3}\triangleq U_{i2}-K_{i1}X_{i2}$. By Theorem 1.7 in \cite{huang2004nonlinear}, the full-information controller $u_i=K_{i1}x_i+K_{i2}d_i+K_{i3}r$ trivially solves the output regulation problem of the $i$th subsystem \eqref{sys:virtual-regulation}, and hence achieves the leader-following coordination goal of the whole multi-agent system in a centralized setup. { Inspired by the separate principle for single linear systems \cite{wonham1979linear}, we follow a similar design and replace the unavailable quantities in the full-information control law by their estimations.}

Since not all agents can directly get access to the reference signal (i.e., the leader), we first construct the following distributed observer for agent $i$ to estimate $r$, and transform the original coordination problem into several decentralized ones:
\begin{align}\label{eq:observer-leader}
 \dot{\eta}_i=S_0\eta_i+ L_0F_0\eta_{vi},
\end{align}
where $\eta_{vi}=\sum_{j=0}^N a_{ij}(t)(\eta_i-\eta_j)$, $\eta_0=r$, $i=1,\dots,N$, and $L_0$ is a constant matrix to be designed.  Letting $\bar \eta_i\triangleq\eta_i-r$ and denoting  $\bar \eta=\mbox{col}(\bar \eta_1,\dots,\bar \eta_N)$ gives
\begin{align}\label{sys:observer-leader}
\dot{\bar \eta}=[I_N \otimes S_0+H_{\sigma(t)}\otimes (L_0F_0)]\bar \eta.
\end{align}

The following lemma shows the effectiveness of this distributed observer.
\begin{lemma}\label{lem:observer-leader}
  Under Assumption \ref{ass:graph}, there exists a constant matrix $L_0$ such that the system \eqref{sys:observer-leader} is uniformly exponentially stable in the sense of $||\bar \eta||\leq c_0  e^{-\lambda_0 t}$ for some positive constants $c_0$ and $\lambda_0$.
\end{lemma}
\noindent{\bf Proof}.  For this purpose, we only have to determine an $L_0$ such that, for each $i$, there exist two constants $\bar c_{0i}$ and $\bar \lambda_{0i}$ such that $||\eta_i-v||\leq c_{0i}  e^{-\lambda_{0i} t}$.

Note that $H_{\sigma(t)}$ is positive definite and constant during each interval $[t_i,t_{i+1})$  under Assumption \ref{ass:graph}. We first consider this problem in each interval. Assume $\sigma(t)=p$ for $t\in[t_i, t_{i+1})$, there exists a unitary matrix $U_p$ such that $\Lambda_p\triangleq U_p^TH_p U_p=\text{diag}\{\lambda_1^p,\dots,\lambda_N^p\}$. Let $\hat \eta=(U_p^T\otimes I_N)\bar \eta$, then,
\begin{align*}
  \dot{\hat \eta}=(I_N\otimes S_0+\Lambda_p\otimes L_0F_0)\hat \eta
\end{align*}
that is, $ \dot{\hat \eta}_i=(S_0+\lambda_i^p L_0F_0)\hat \eta_i$ for $i=1,\dots,N$, where $\lambda_i^p>0$ for $i=1,...,N$ are the eigenvalues of $H_p$ ($p\in \mathcal{P}$). Since $(S_0,\, F_0)$ are detectable, there exists \cite{boyd1994linear} a positive definite symmetric matrix $P$ satisfying
\begin{equation}\label{eq:lyapunov}
PS_0+S_0^TP-2F_0^TF_0<0.
\end{equation}
Note that the minimum eigenvalue of $H_p$ for all $p$ is well-defined. Denoting it as $\bar \lambda>0$ and taking  $L_0=-\mu^* P^{-1}F_0^T$ with $\mu^*\triangleq\max\{\frac{1}{\bar \lambda},1\}$ gives
\begin{align*}
  &(S_0+\lambda_i^p L_0F_0)^T P+ P(S_0+\lambda_i^p L_0F_0)\\
  &=S_0^TP+PS_0-2\mu^*\lambda_i^p F_0^TF_0\\
  &=S_0^TP+PS_0-2F_0^TF_0-2(\mu^*\lambda_i^p-1) F_0^TF_0\\
  &\leq S_0^TP+PS_0-2F^T_0F_0
\end{align*}
Since $S_0^TP+PS_0-2F_0^TF_0$ is negative definite, under Assumption \ref{ass:graph}, there exists a positive constant $c$, such that
\begin{equation}\label{eq:lmi}
(S_0+\lambda_i^p L_0F_0)^T P+ P(S_0+\lambda_i^p L_0F_0)\leq -cP
\end{equation}

By letting $V_{\eta}=\sum_{i=1}^N{\hat \eta}_i^T P {\hat \eta}_i$, we can derive
$\dot{V}_{\eta}\leq -c V_{\eta}$. Recalling the dwell-time assumption this inequality holds for all $t$. Note that $\bar \eta^T \bar \eta=\hat \eta^T \hat \eta$, it follows
\begin{align*}
  ||\eta_i-v||^2\leq\hat \eta^T\hat \eta \leq \lambda_{min}(P)^{-1}{V}_{\eta}(t)< \lambda_{min}(P)^{-1}V_{\eta}(0)e^{-ct}
\end{align*}
The conclusion is readily obtained.
\hfill\rule{4pt}{8pt}

\begin{remark}
  Although $\eta_0=r$ appears in \eqref{eq:observer-leader}, $y_0=F_0\eta_0$ will suffice this design. When $y_0=r$ (i.e., $F_0=I_p$), it means that the state of the leader can be directly obtained when some agent is connected to it. This circumstance has been partly considered in \cite{su2012scl}.  We extend these results using only output measurements of the leader to deal with the cases when only partial states are available. Similar control laws were proposed in \cite{hong2006tracking} when the followers are all integrators, while here we consider general linear agents and also local disturbances.
\end{remark}

After building distributed observers for those followers, it is natural to replace $r$ by its estimation $\eta_i$.  The following lemma guarantees the validity of this substitution and shows how it transforms the coordination problem of these multi-agent systems into $N$ decentralized estimation and regulation sub-tasks.

\begin{lemma}\label{lem:decentralized}
 Under Assumptions \ref{ass:graph} and \ref{ass:regulator-equation}, $u_i=K_{i1}x_i+K_{i2}d_i+K_{i3}\eta_i, \,\dot{\eta}_i=S_0\eta_i+ L_0F_0\eta_{vi}$ will solve the output regulation problem of system \eqref{sys:virtual-regulation}, and hence the leader-following coordination problem of this multi-agent system with the selected $L_0$, where $K_{i1}$, $K_{i2}$, $K_{i3}$ are matrices defined in the centralized case.
\end{lemma}
\noindent{\bf Proof}.   Under Assumption \ref{ass:regulator-equation}, letting $\bar x_i=x_i-X_{i1}d_i-X_{i2}r$ gives
\begin{align}\label{sys:closed1-lem2}
\begin{split}
  \dot{\bar x}_i&=(A_i+B_iK_{i1})\bar x_i+B_iK_{i3}\bar \eta_i\\
  \dot{\bar\eta}_i&=S_0\bar \eta_i+L_0F_0\eta_{vi}\\
  e_i&=(C_i+D_iK_{i1})\bar x_i+D_iK_{i3}\bar \eta_i
  \end{split}
\end{align}
or in compact form
\begin{align}\label{sys:closed2-lem2}
\begin{split}
  \dot{\bar x}&=\bar A\bar x+ \bar B \bar \eta\\
  \dot{\bar\eta}&=(I_N\otimes S_0+H_{\sigma(t)}\otimes L_0F_0)\bar \eta\\
  e&=\bar C\bar x+\bar D\bar \eta
  \end{split}
\end{align}
where $\bar x=\mbox{col}(x_1,\dots,x_N)$ and
\begin{align*}\bar A=\mbox{block}\,diag\{A_1+B_1K_{11},\dots,A_N+B_NK_{N1}\},\quad \bar B=\mbox{block}\,diag\{B_1K_{13},\dots,B_NK_{N3}\},\\
\bar C=\mbox{block}\,diag\{C_1+D_1K_{11},\dots, C_N+D_NK_{N1}\},\;\bar D=\mbox{block}\,diag\{D_1K_{13},\dots,D_NK_{N3}\}.\end{align*}
Since $\bar A$ is Hurwitz, by Lemmas \ref{lem:time-varying} and \ref{lem:observer-leader}, $e_i=y_i-y_0$ will converge to zero as $t\to \infty$.
\hfill\rule{4pt}{8pt}

Next, we aim to solve those decentralized reference tracking and disturbance rejection problems. As having been pointed out before that the disturbances are often unmeasurable, the control law in Lemma \ref{lem:decentralized} is not implementable. To tackle this problem, we employ the disturbance observer based (DOB) approach which has been well-studied in its centralized version by many authors \cite{nakao1987robust,li2014disturbance}.

Basically, we seek to reconstruct an estimation of the disturbances that affect the tracking performance, and then use these estimations to achieve disturbance rejection.  Hence, the following assumption comes naturally.
\begin{assumption}\label{ass:observable}
For each $i=1,\dots,N$, the pair $\left([C_i, \,\mathbf{0}^{l\times q_i}],\begin{bmatrix}
  A_i&E_i\\
  0_{n_i}&S_i
\end{bmatrix}\right)$ is detectable.
\end{assumption}

We first consider the output feedback cases, where only $y_i$ is available through measurement, and propose a composite control law in the following form.
\begin{align}\label{ctr:law1}
     \dot{\xi}_i&=A_i\xi_i+E_i\zeta_i+ B_iu_i-L_{i1}(y_i-\hat y_i) \nonumber\\
     \dot{\zeta}_i&= S_i\zeta_i-L_{i2}(y_i-\hat y_i)\nonumber\\
     \dot{\eta}_i&=S_0\eta_i+ L_0F_0\eta_{vi}\nonumber\\
      u_i&=K_{i1}\xi_i+K_{i2}\zeta_i+K_{i3}\eta_i, \quad i=1,\dots, N
\end{align}
where $\hat y_i=C_i\xi_i+D_iu_i$, $i=1,\dots,N$, and $L_{i1}$, $L_{i2}$ are constant matrices such that $$A_{ci}\triangleq \begin{bmatrix}
  A_i+L_{i1}C_i&E_i\\
  L_{i2}C_i&S_i
\end{bmatrix}
$$ is Hurwitz.

It is time to  give our first main theorem.
\begin{theorem}\label{thm1}
  Under Assumptions \ref{ass:graph}-\ref{ass:observable}, the leader-following coordination problem of the multi-agent system composed of \eqref{follower}, \eqref{disturbance}, and \eqref{leader} can be solved by the control law \eqref{ctr:law1}.
\end{theorem}
\noindent{\bf Proof}.  Let $\hat x_i=\xi_i-x_i$, $\hat d_i=\zeta_i-d_i$, and $\bar x_i=x_i-X_{i1}d_i-X_{i2}r$. Under the control law \eqref{ctr:law1}, the closed-loop system of agent $i$ is of the form:
\begin{align}\label{sys:closed1-thm1}
\begin{split}
  \dot{\bar x}_i&=(A_i+B_iK_{i1})\bar x_i+B_iK_{i1}\hat x_i+ B_iK_{i2}\hat d_i+ \bar B_iK_{i3}\bar \eta_i \\
  \dot{\hat x}_i&=(A_i+L_{i1}C_i)\hat x_i+E_i\hat d_i\\
  \dot{\hat d}_i&=L_{i2}C_i\bar x_i+S_i\hat d_i\\
  \dot{\bar \eta}_i&=S_0\bar \eta_i+ L_0F_0\eta_{vi} \\
  e_i&=(C_i+D_iK_{i1})\bar x_i+D_iK_{i1}{\hat x}_i+D_iK_{i2}\hat d_i+D_iK_{i3}\bar \eta_i
  \end{split}
\end{align}
Let $\bar x=\mbox{col}(x_1,\dots,x_N)$, $\bar v=\mbox{col}\{\hat x_1,\,\dots,\,\hat x_N,\,\hat d_1,\,\dots,\,\hat d_N,\,\bar \eta_1,\,\dots,\,\bar \eta_N\}$, and the whole multi-agent system can be put into a compact form as
\begin{align}\label{sys:closed2-thm1}
\begin{split}
  \dot{\bar x}&=\bar A \bar x+ \bar B \bar v\\
  \dot{\bar v}&=\bar S \bar v\\
  e&=\bar C\bar x+\bar F\bar v
  \end{split}
\end{align}
where $\bar A=\mbox{block}\,diag\{A_1+B_1K_{11},\dots,A_N+B_NK_{N1}\},$\,$\bar S=\mbox{block}\,diag\{A_{c1},\dots,A_{cN},I_N\otimes S+H_{\sigma(t)}\otimes L_0F_0\},$\,$\bar B=[\bar  B_1, \bar B_2,\bar B_3]$,\,$\bar C=\mbox{block}\,diag\{C_1+D_1K_{11},\,\dots,\,C_N+D_NK_{N1}\},$\,$\bar F= [\bar F_1,\, \bar F_2,\,\bar F_3]$ and $\bar B_k=\mbox{block}\,diag\{B_1K_{1k},\,\dots,\,B_NK_{Nk}\}$, $\bar F_k=\mbox{block}\, diag\{D_1K_{1k}, \,\dots ,\,D_NK_{Nk}\}$ ($k=1,2,3$). Since $\bar A$ and $A_{ci}$ are Hurwitz, by Lemmas \ref{lem:time-varying} and \ref{lem:observer-leader}, $e_i=y_i-y_0$ will converge to zero as $t\to \infty$. The proof is thus completed.
\hfill\rule{4pt}{8pt}

\begin{remark}
When $E_i=0$ for all agents, it reduces to the well-studied leader-following consensus problem considering only reference tracking problem. Then, the relevant results in \cite{olfati2007consensus} and \cite{ren2008distributed} are actually special cases of this theorem for integrators. Even when the leader has a general linear dynamics as in \cite{ni2010leader}, we consider both reference tracking and local disturbance rejection problems under switching topologies.
\end{remark}
\begin{remark}\label{rem:dobc}
When $N=1$, this problem becomes a centralized reference tracking and disturbance rejection problem. While most of existing DOB results focus on rejecting those disturbances \cite{chen2000nonlinear, li2014disturbance}, our design also incorporates the reference tracking aspects. Those DOB controllers can admit not only bounded disturbances (e.g., constants and sinusoidal signals) but also unbounded disturbances (e.g., ramping signals and polynomials) under switching topologies.
\end{remark}

In many circumstances, the state $x_i$ may be available for us by direct measurement or other treatments. Thus, there exists a certain degree of redundancy in the controller \eqref{ctr:law1}, which still produces the estimations of $x_i$. To remove such redundancies and save our controller's order, we propose a reduced-order DOB control law to facilitate our design.

For this multi-agent system, the reduce-order disturbance observer is given as follows.
\begin{align}\label{ctr:law2}
      \dot{\zeta}_i&=(S_i+L_iE_i)\zeta_i+(L_iA_i-S_iL_i-L_iE_iL_i)x_i+L_iB_iu_i\nonumber\\
      \dot{\eta}_i&=S_0\eta_i+ L_0F_0\eta_{vi}\nonumber\\
      u_i&=K_{i1}x_i+K_{i2}(\zeta_i-L_ix_i)+K_{i3}\eta_i, \quad i=1,\dots, N
\end{align}
where $K_{ij}$, $i=1,\dots,N,\,j=1,2,3$, and $L_i$ are gain matrices to be determined later.

With this reduced-order controller, the following theorem can be derived.
\begin{theorem}\label{thm2}
  Under Assumptions \ref{ass:graph}--\ref{ass:observable}, there exist constant matrices $K_{ij}$ and $L_i$, $i=1,\dots, N,\, j=1,\,2,\,3$, such that the leader-following coordination problem of this multi-agent system is solved by the control law \eqref{ctr:law2}.
\end{theorem}
\noindent{\bf Proof}.  The proof is similar with that of Theorem \ref{thm1}.
Let $\hat d_i=\zeta_i-L_ix_i-d_i$, and $\bar x_i=x_i-X_{i1}d_i-X_{i2}r$. Under the control law \eqref{ctr:law2}, the closed-loop system of agent $i$ is with the form of
\begin{align}\label{sys:closed1-thm2}
\begin{split}
  \dot{\bar x}_i&=(A_i+B_iK_{i1})\bar x_i+ B_iK_{i2}\hat d_i+ \bar B_iK_{i3}\bar \eta_i \\
  \dot{\hat d}_i&=(S_i+L_iE_i)\hat d_i\\
  \dot{\bar \eta}_i&=S_0\bar \eta_i+ L_0F_0\eta_{vi} \\
  e_i&=(C_i+D_iK_{i1})\bar x_i+D_iK_{i2}\hat d_i+D_iK_{i3}\bar \eta_i
  \end{split}
\end{align}
Let $\bar x=\mbox{col}(x_1,\dots,x_N)$, $\bar v=\mbox{col}\{\hat d_1,\,\dots,\,\hat d_N,\,\bar \eta_1,\,\dots,\,\bar \eta_N\}$. By some mathematical manipulations, the whole multi-agent system can be put into a compact form as
\begin{align}\label{sys:closed2-thm2}
\begin{split}
  \dot{\bar x}&=\bar A \bar x+ \bar B \bar v\\
  \dot{\bar v}&=\bar S \bar v\\
  e&=\bar C\bar x+\bar F\bar v
  \end{split}
\end{align}
where
\begin{align*}
\bar A&=\mbox{block}\,diag\{A_1+B_1K_{11},\dots,A_N+B_NK_{N1}\}\\
\bar B&=\mbox{block}\,{diag}\{B_1K_{13},\,\dots,\,B_NK_{N3}\},\\
\bar C&=\mbox{block}\,diag\{C_1+D_1K_{11},\,\dots,\,C_N+D_NK_{N1}\},\\
\bar F&=\mbox{block}\,{diag}\{D_1K_{13},\,\dots,\,D_NK_{N3}\},\\
\bar S&=\mbox{block}\,diag\{S_1+L_1E_1,\dots,S_N+L_NE_N,I_N\otimes S+H_{\sigma(t)}\otimes L_0F_0\}.
\end{align*}

According to Lemma \ref{lem:time-varying} and \ref{lem:observer-leader},we only have to find proper $K_{ij}$ and $L_i$ such that $\bar A$ and $S_i+L_iE_i$ are all Hurwitz. Then, by similar arguments as that in Theorem \ref{thm1}, $e_i=y_i-y_0$ will converge to zero as $t\to \infty$.  In fact, such gain matrices indeed exist. Take $K_{ij}$ as defined in Theorem \ref{thm1}, and the detectability of $(E_i,\,S_i)$ will suffice the selection of $L_i$, which is obvious by PBH-test under Assumption \ref{ass:observable}.  Thus the conclusion follows readily.
\hfill\rule{4pt}{8pt}

\begin{remark}\label{rem:disturbance}
  As having been pointed before, unlike in existing cooperative output regulation result \cite{su2012scl}, the disturbances are locally modeled by different autonomous systems from that of the global reference. A similar setup has been used in \cite{wang2014semi} and \cite{tang2015auto}. This separate modeling method results in two dissimilar treatments, distributed observers for the global reference and DOB approach for the local disturbances, which may enhance the effectiveness of our design and bring a better performance. Also, even for the case when $N=1$, we proposed different reduced-order disturbance observers from that in \cite{li2014disturbance} to solve this problem without using the derivative of the plant's states.
\end{remark}

\section{Simulations}

\begin{figure}
  \centering
  \subfigure[The graph $\mathcal{G}_1$]
    {\centering
\begin{tikzpicture}[shorten >=1pt, node distance=1.2 cm, >=stealth',
every state/.style ={circle, minimum width=0.2cm, minimum height=0.2cm}, auto]
\node[align=center,state](node0) {0};
\node[align=center,state](node1)[right of=node0]{1};
\node[align=center,state](node2)[right of=node1]{2};
\node[align=center,state](node3)[right of=node2]{3};
\path[->]   (node0) edge (node1)
            (node1) edge [bend right] (node2)
            (node2) edge [bend right] (node1)
            (node0) edge [bend left]  (node3)
            ;
\end{tikzpicture}
  }\quad
  \subfigure[The graph $\mathcal{G}_2$]
   {\centering
\begin{tikzpicture}[shorten >=1pt, node distance=1.2 cm, >=stealth',
every state/.style ={circle, minimum width=0.2cm, minimum height=0.2cm}, auto]
\node[align=center,state](node0) {0};
\node[align=center,state](node1)[right of=node0]{1};
\node[align=center,state](node2)[right of=node1]{2};
\node[align=center,state](node3)[right of=node2]{3};
\path[->]   (node0) edge (node1)
            (node0) edge [bend left](node2)
            (node2) edge [bend right] (node3)
            (node3) edge [bend right] (node2)
            ;
\end{tikzpicture}
    }
\caption{The communication graphs.}\label{fig:graph}
\end{figure}
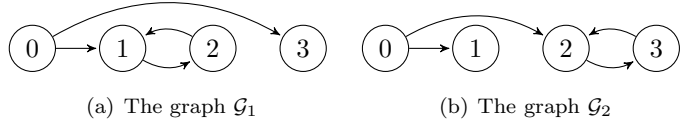

As an example, we consider the coordination problem for a multi-agent system consisting of three followers and one leader. The follower agents are the mass-damper-spring systems with unit
mass described by:
\begin{align*}
  \ddot{y}_i+g_i\dot{y}_i+f_iy_i=u_i+E_id_i,\quad i=1,2,3
\end{align*}
where $d_i$ is the local disturbance. Those disturbances are modeled by $S_1=[0,1;0,0], E_1=[1,0]$; $S_2=0, E_2=1$; $S_3=[0,1;-1,0],\,E_3=[1,0]$. The leader is specified by a harmonic oscillator:
$\dot{r}_1=r_2,\,\dot{r}_2=-r_1, \, y_0=r_1$. We assume here the interconnection topology is switching between graph $\mathcal{G}_1$ and $\mathcal{G}_2$ described by Fig.~\ref{fig:graph}. The switchings are periodically carried out in the following order $\{\mathcal{G}_1, \mathcal{G}_2, \mathcal{G}_1, \mathcal{G}_2, \cdots\}$ with switching period $t=5$.

Letting $x_{i1}=y_i,\, x_{i2}=\dot{y}_i$, then,
$$\dot{x}_{i1}=x_{i2},\quad \dot{x}_{i2}=-f_ix_{i1}-g_ix_{i2}+u_i+E_id_i, \quad y_i=x_{i1}.$$
Apparently, the coordination laws in \cite{hong2006tracking} and \cite{ni2010leader} will not work for these agents. Even the discontinuous rule (\cite{bauso2009consensus}) fails to solve this problem because of unbounded disturbances in Agent $1$. Nevertheless, as Assumptions \ref{ass:graph}-\ref{ass:observable} are satisfied, it can be solved by the methods given in last sections.

\begin{figure}
  \centering
  \subfigure[]{
    \includegraphics[width=0.4\textwidth]{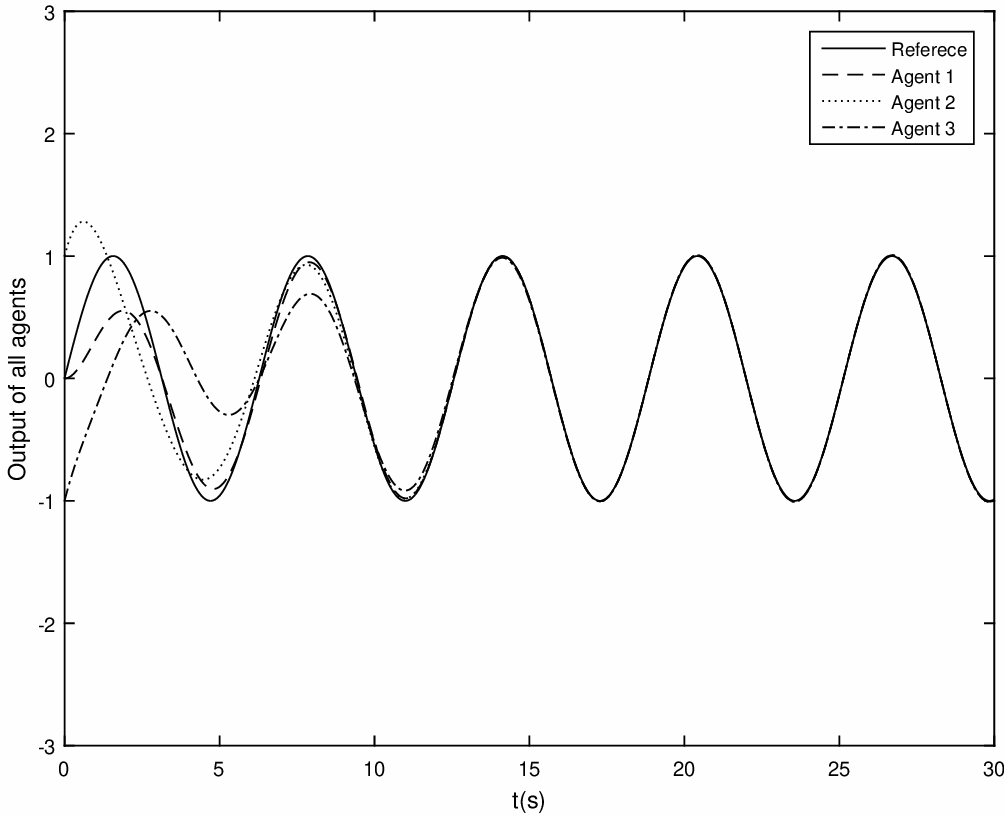}
  }\quad
  \subfigure[]{
    \includegraphics[width=0.4\textwidth]{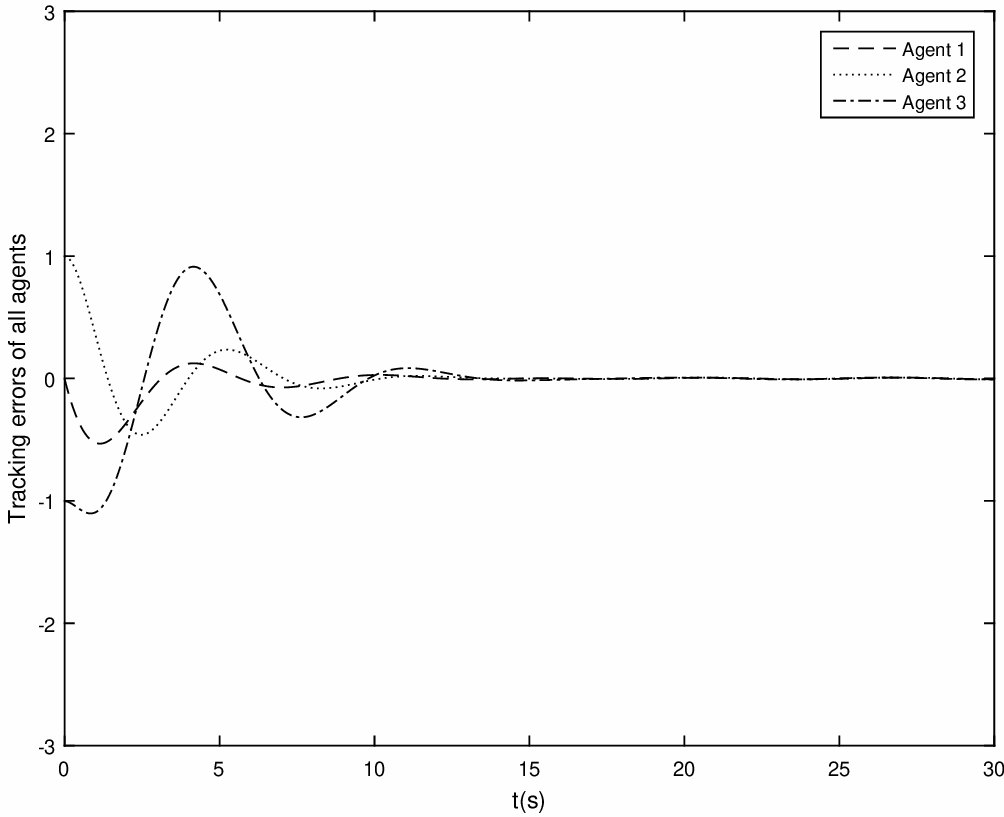}
    }
\caption{Tracking performance with output feedback control law \eqref{ctr:law1}.}\label{fig:output}
\end{figure}

\begin{figure}
  \centering
  \subfigure[]{
    \includegraphics[width=0.4\textwidth]{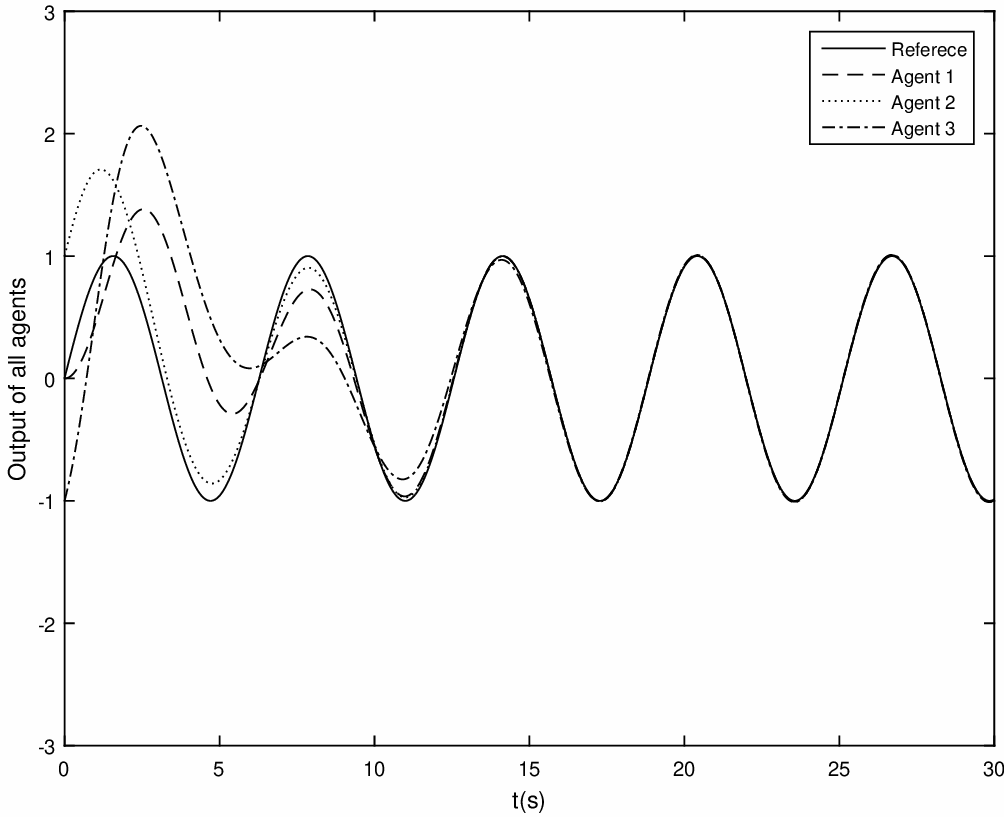}
  }\quad
  \subfigure[]{
    \includegraphics[width=0.4\textwidth]{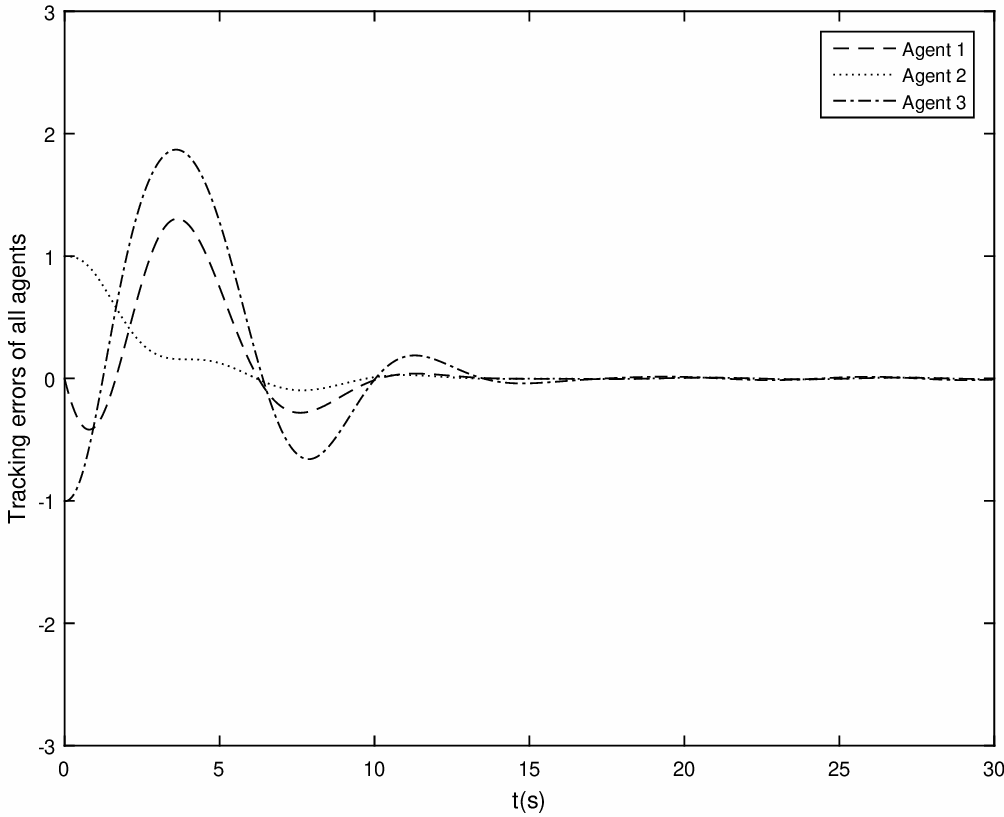}
    }
\caption{Tracking performance with reduced-order DOB feedback control law \eqref{ctr:law2}.}\label{fig:reduced}
\end{figure}

For simulations, the system parameters are taken as $f_1=1, g_1=1, f_2=0, g_2=1$ and $f_3=1,g_3=0$. By solving the regulator equations in Assumption \ref{ass:regulator-equation} and also the Lyapunov inequality \eqref{eq:lyapunov}, we choose proper gain matrices for the controllers as showed in Table~\ref{gain-table}. While the initials for the plant is generated between $[-1,\,1]^2$, the initials for controllers are set at their origins. The simulation results using full-order and reduced-order disturbance observer based control are showed in Fig.~\ref{fig:output} and Fig.~\ref{fig:reduced}, respectively.  We list the outputs of agents at serval time points in Table~\ref{output-table}. It can be found that after $21s$, the agents can track the leader and reject those disturbances with an error tolerance $2\times 10^{-3}$.

\begin{table}
\centering
\small
 \begin{tabular}{c|c|c}
  \hline
  $\quad$ & {\footnotesize Simulation 1} & {\footnotesize Simulation 2}\\
  \hline
  \multirow{2}{*}{\footnotesize Agent 1} &\multicolumn{2}{c}{$K_{11}=[0, 0]\quad K_{12}=[-1,0]\quad K_{13}=[0,1]\quad L_0=[-1.3522;-0.4142]$} \\
  \cline{2-3}
  &$L_{11}=[-1.9813; -1.4628] \quad L_{12}=[-3.1445; -1.0000]$ & $L_1=[0,-1;0,-1]$\\
  \hline
  \multirow{2}{*}{\footnotesize Agent 2} &\multicolumn{2}{c}{$K_{21}=[-1 , 0]\quad K_{22}=-1\quad K_{23}=[0,1]\quad L_0=[-1.3522;-0.4142]$} \\
\cline{2-3}
  &$L_{21}=[-1.7321; -1.0000]\quad L_{22}=-1.0000$ & $L_2=[0, -1]$\\
  \hline
  \multirow{2}{*}{\footnotesize Agent 3} &\multicolumn{2}{c}{$K_{31}=[0,-1]\quad K_{32}=[-1 ,0]\quad K_{33}=[0 ,1], L_0=[-1.3522;-0.4142]$} \\
\cline{2-3}
  &$L_{31}=[-2.1607;-1.8343] \quad L_{32}=[-0.8860; 1.1023]$& $L_{3}=[0,-1; 0,0]$\\
  \hline
\end{tabular}
\caption{Gain matrices for each agent in simulations.}\label{gain-table}
\end{table}

\begin{table}
\centering
\small
 \begin{tabular}{c|c|cccccccc}
  \hline
   \multicolumn{2}{c|}{\footnotesize Sampling time}& t=6s & t=12s & t=15s & t=18s & t=21s & t=24s& t=27s&t=30s \\
  \hline
  \multicolumn{2}{c|}{\footnotesize Reference output}&  -0.2794   &-0.5366   & 0.6503  & -0.7510   & 0.8367   &-0.9056    &0.9564   &-0.9880 \\
  \hline
  \multirow{3}{*}{\footnotesize Simulation 1} & Agent 1 & -0.8507   &-0.6008  & 0.6295   & -0.7574    &0.8347  & -0.9064   & 0.9562   &-0.9882\\
  \cline{2-10}
  &Agent 2&  -0.6961   &-0.5803 &  0.6395   &-0.7532    &0.8365   &-0.9057  &  0.9565 &  -0.9881\\
  \cline{2-10}
  &Agent 3&  -0.3552  & -0.5100   &0.6400   &-0.7494   & 0.8364 &  -0.9058  &  0.9565 &  -0.9882\\
  \hline
  \multirow{3}{*}{\footnotesize Simulation 2} & Agent 1 &    -0.3054 &  -0.6109  &  0.6292   &-0.7574   & 0.8347   &-0.9064  &  0.9562   &-0.9882\\
  \cline{2-10}
  &Agent 2&     -0.6476   &-0.5807  &0.6395   &-0.7532   & 0.8365 &  -0.9057  &  0.9565   &-0.9881\\
  \cline{2-10}
  &Agent 3&    -0.1382  & -0.4240   &0.6205   &-0.7466   & 0.8363  & -0.9058  &  0.9566   &-0.9882 \\
  \hline
\end{tabular}\caption{Outputs of each agent in simulations.}\label{output-table}
\end{table}

\section{Conclusions}
A leader-following coordination problem was solved for a class of heterogeneous multi-agent systems subject to local disturbances under switching topologies. By devising a distributed observer, this problem was transformed into several decentralized estimation and regulation sub-tasks, and eventually solved by two disturbance observer based control laws. Our future work will include nonlinear cases and with more general graphs.

\bibliographystyle{ieeetran}
\bibliography{mybibfile}
\nocite{*}

\end{document}